\documentclass[final]{autart}

\usepackage{amssymb}
\usepackage{makeidx}
\usepackage{multicol}
\usepackage[bottom]{footmisc}
\usepackage{graphicx}
\usepackage{color}
\usepackage{natbib}
\usepackage{epstopdf}
\usepackage{url}

\newcommand{\eq}{\begin{equation}}
\newcommand{\eeq}{\end{equation}}
\newcommand{\eqn}{\begin{eqnarray}}
\newcommand{\eeqn}{\end{eqnarray}}

\newcommand{\bsea}{\begin{subeqnarray}}
\newcommand{\esea}{\end{subeqnarray}}
\newcommand{\nn}{\nonumber}

\newcommand{\tr}{\mathop{\rm tr}}  

\newcommand{\Ac}{ \mathcal{A}}
\newcommand{\Bc}{ \mathcal{B}}

\newcommand{\Dc}{ \mathcal{D}}

\newcommand{\Mc}{ \mathcal{M}}

\newcommand{\Pc}{ \mathcal{P}}
\newcommand{\Qc}{ \mathcal{Q}}

\newcommand{\Sc}{ \mathcal{S}}
\newcommand{\Tc}{ \mathcal{T}}

\newcommand{\Vc}{ \mathcal{V}}

\newcommand{\Cs}{ \mathbb{C}}

\newcommand{\Es}{ \mathbb{E}}

\newcommand{\Ns}{ \mathbb{N}}

\newcommand{\Rs}{ \mathbb{R}}

\newcommand{\Ts}{ \mathbb{T}}

\newcommand{\Zs}{ \mathbb{Z}}

\def\qed{\hfill \vrule height 7pt width 7pt depth 0pt \smallskip}

\begin{document}

\begin{frontmatter}
\title{An Interpretation of the Dual Problem\\ of the THREE-like Approaches \thanksref{footnoteinfo}}
\thanks[footnoteinfo]{This work has been partially supported by the FIRB project ``Learning
meets time'' (RBFR12M3AC) funded by MIUR.}


\author[Padova]{Mattia Zorzi}\ead{zorzimat@dei.unipd.it}

\address[Padova]{Dipartimento di Ingegneria dell'Informazione, Universit\`a degli studi di
Padova, via Gradenigo 6/B, 35131 Padova, Italy}

\begin{keyword}
Spectral Estimation; Prediction Error Identification Method; Convex Optimization; Divergence Family.
\end{keyword}

\begin{abstract}
Spectral estimation can be preformed using the so called THREE-like approach. Such method leads to a convex optimization problem whose solution is characterized through its dual problem. In this paper, we show that the dual problem can be seen as a new parametric spectral estimation problem. This interpretation implies that the THREE-like solution is optimal in terms of closeness to the correlogram over a certain parametric class of spectral densities, enriching in this way its meaningfulness.
\end{abstract}
\end{frontmatter}

\section{Introduction}
In science and engineering it is often required to approximate phenomena through simple models. The quality of this approximation heavily depends on its intended use. For instance, if the purpose is prediction, then the prediction error identification method (PEM)
provides the optimal model, see \cite{LJUNG_IDENTIFICATIONTHEORYUSER}, \cite{SODERSTROM_STOICA} and \cite{LINDQUIST_PICCI}. Therefore, the more interpretations the solution of an identification method admits, the more importance such a solution gains.

In this paper, we deal with a spectral estimation problem called THREE-like method, see the former work in \cite{BYRNES_GEORGIOU_LINDQUIST_THREE}.
Such approach exploits the output covariance matrix of a filter to extract information on the spectral density of the input process. More precisely, the class of input spectral densities matching the output covariance matrix is considered and a spectrum approximation problem, which chooses an estimate of the input spectral density in this class, is then employed. Such choice is the spectral density minimizing a divergence index with respect to an {\em a priori} spectral density.   
 The strength of this approach is the flexibility in choosing the filter, the divergence index and the {\em a priori} spectral density. In particular, choosing properly them, it is possible to recover the maximum entropy estimator \citep{BURG_MAX_ENTROPY_SPECTRAL_ANALYS} used to estimate autoregressive (AR) models.

In this paper, 
we show that the dual problem of the spectrum approximation problem can be seen as a new parametric spectral estimation problem wherein the optimal solution is the closest spectral density to the {\em correlogram}, according to a suitable weighted divergence index, and belonging to a certain parametric class. Therefore, this new interpretation enriches the meaningfulness of the optimal solution. Moreover, two specific THREE-like solutions can be also interpreted as solution to PEM. This interpretation has been observed formerly in \cite{LINDQUIST_PREDICTION_ERROR_2007}  and \cite{BLOMQVIST_WAHLBERG_2007} for a particular setting.

The outline of the paper follows. In section \ref{sec_review_three} we review the THREE-like approaches. In Section \ref{sec_weighted_beta} we define two types of weighted divergence index needed to introduce, in Section \ref{sec_dual_interpret}, the interpretation of the dual problem. Finally, in Section \ref{sec_PEM} we show the link between THREE and PEM.

Throughout the paper we use the following notation.
$\Qc_n$ denotes the vector space of $n\times n$ symmetric matrices, and $\Ns_+=\Ns\setminus \{0\}$.
We drop the dependence of the domain for functions which are defined over the unit circle, for instance for spectral densities. If a function $\Phi$ is positive (semi)definite on the unit circle we write $\Phi\succ 0$ ($\Phi\succeq 0$). The shorthand notation $\int \Phi$ means the integration of $\Phi$ over the unit circle with respect to the normalized {\em Lebesgue} measure.

\section{A review of the THREE-like Approaches}\label{sec_review_three}
Let $y=\{y(t)\}_{t\in\Zs}$ be a $\Rs^m$-valued zero-mean stationary purely nondeterministic Gaussian process.
Such a process is completely characterized by its spectral density denoted by $\Phi$. Recall that, $\Phi$ is a $m\times m$ Hermitian-valued positive semidefinite function defined over the unit circle. We assume that $\Phi$ belongs to the set $\Sc_m^+$ of spectral densities which are bounded and coercive on the unit circle, that is there exist two scalar constants $k_1,k_2>0$ such that $k_1 I\preceq \Phi\preceq k_2 I$.

A THREE-like approach is a procedure to estimate $\Phi$ from a finite length sequence $\mathrm{y}=\{\mathrm{y}(t)\}_{t=1}^N$ extracted from a realization of $y$.
It can be sketched as follows:
\begin{itemize}
\item Fix a filter $G(z)=(zI-A)^{-1}B$ with $A\in\Rs^{n\times n}$ strictly stable, $B\in\Rs^{n\times m}$, $n> m$, and such that the pair $(A,B)$ is reachable. In this way, $\Sigma:=\Es[x(t)x(t)^T]\succ 0$ where $x=\{x(t)\}_{t\in\Zs}$ is the zero-mean stationary Gaussian process at the output of the filter $G$ when fed by $y$. Then, compute an estimate $\hat \Sigma$ of $\Sigma$ from $\mathrm{y}$.
\item Fix an {\em priori} spectral density (i.e. prior) $\Psi\in\Sc_m^+$ for $y$ with bounded {\em McMillan} degree. More precisely, $\Psi$ is designed using some given partial information 
on $y$ (e.g. its zeroth moment) or using given laws describing theoretical features of $y$. Note that, the specification of the {\em prior} $\Psi$ is not strictly necessary: if no a {\em priori} information is available, we set $\Psi=I$ corresponding to white Gaussian noise with variance equal to the identity (WGN), i.e. the most unpredictable process. 
\item The estimate of $\Phi$ is given by solving the following spectrum approximation problem
\eqn \label{approx_pb}&&  \Phi^\circ=\mathrm{argmin}_{\Phi\in\Sc_m^+} \Sc(\Phi\|\Psi)\nn\\
    && \hspace{1cm} \hbox{s.t. } \int G \Phi G^*=\hat \Sigma\eeqn
where $\Sc$ is a pseudo-distance (or divergence index) between two spectral densities in $\Sc_m^+$, that is
$\Sc(\Phi\|\Psi)\succeq 0$ for any $\Phi,\Psi\in\Sc_m^+$ and equality holds if and only if $\Phi=\Psi$. \end{itemize}

Thus, $\Phi^\circ$ is the closest spectral density to $\Psi$, according to $\Sc$, matching the estimated output covariance matrix $\hat\Sigma$. 
The THREE-like approach is a generalization of the maximum entropy 
estimator used for AR modeling. Indeed, if we choose $G$ as a bank of $n$ delays, $\Psi=I$ and $\Sc$ the {\em Kullback-Leibler} divergence (see below), then $\Phi^\circ$ is the spectral density matching the first $n$ estimated covariance lags of $y$ and maximizing the entropy rate 	\citep{COVER_THOMAS}.   
It is also worth noting that Problem (\ref{approx_pb})
can be seen as a rigorous recasting of the beamspace technique used to determine the direction of arrival of narrow-band emitter signals impinging on an antenna array \citep{BEAMSPACE_APPROACH}. In that case, $y$ is the observation vector, $x$ is the beamspace data vector and $G$ is the beamforming matrix.

The filter $G$ is designed by the user to represent known dynamics or to post-process the data $\mathrm{y}$ . In the latter case, this freedom represents a powerful tool to perform high resolution spectral analysis. More precisely, a higher resolution can be attained by selecting the poles of $G$ in the proximity of the unit circle,  with arguments in the range of interest \citep{BYRNES_GEORGIOU_LINDQUIST_THREE, GEORGIOU_STATECOV}. Finally, from $G$ it is also possible to characterize the uncertainty set of $\Phi^\circ$ about the actual spectral density according to a suitable notion of distance \citep{UNCERTAINTY_KARLSSON_2012}.

Regarding the computation of $\hat \Sigma$, we consider the procedure in \cite{ON_THE_ESTIMATION_ZORZI_2012} which is based on the characterization of $\Sigma$ in terms of $G$ and the covariance lags sequence of $y$. It has been shown that $\hat \Sigma \succeq 0$ (and $\hat \Sigma\succ 0$ in all simulations) and such that
\eq \label{cond_correlogramma}\int G \Omega G^* =\hat \Sigma.\eeq
where $\Omega$ is the {\em biased correlogram}, possibly truncated with the {\em Baclman-Tukey method} \citep{STOICA_MOSES_SPECTRALANALYSIS}. Note that, $\Omega\succ 0$ with high probability and it represents a nonparametric spectral estimator of $y$. It is worth noting that $\Sigma$ can be estimated setting a convex optimization problem, see \cite{OTTSIGMA}, however, condition (\ref{cond_correlogramma}) does not hold with the {\em correlogram} but for some unknown spectral density.

The choice of the divergence index $\Sc$ has important implications in the solution $\Phi^\circ$, for instance it fixes the upper bound on the {\em McMillan} degree 
of $\Phi^\circ$. The divergence indexes proposed in the literature can be classified in three different divergence families: the Alpha, Beta and Tau divergence family. In the next sections we review the family of THREE-like solutions with these divergence families. 

\subsection{Solution with the Alpha Divergence Family}
In \cite{ALPHA}, it has been considered the Alpha divergence family
\eqn && \label{def_alpha_divergence}\Sc_A^{(\alpha)}(\Phi\|\Psi)=\tr\int [\frac{1}{\alpha(\alpha-1)}\Phi^\alpha\Psi^{1-\alpha}-\frac{1}{\alpha-1}\Phi\nn\\
&& \hspace{1cm}+\frac{1}{\alpha}\Psi],\;\; \alpha\in\Rs \setminus \{0,1\}.\eeqn
For $\alpha=0$ and $\alpha=1$, it is defined by continuity
\eqn &&\lim_{\alpha\rightarrow 0}\Sc_A^{(\alpha)}(\Phi\|\Psi)=\Sc_{KL}(\Psi\|\Phi) \nn\\
&&\lim_{\alpha\rightarrow 1}\Sc_A^{(\alpha)}(\Phi\|\Psi)=\Sc_{KL}(\Phi\|\Psi)\eeqn
where $\Sc_{KL}$ is the {\em Kullback-Leibler} divergence
\eq \Sc_{KL}(\Phi\|\Psi)=\tr\int [\Phi(\log\Phi-\log\Psi)-\Phi+\Psi].\eeq
By considering the scalar case, $m=1$, and the parametrized family $\Sc_A^{(1-\nu^{-1})}(\Phi\|\Psi)$ with $\nu\in\Ns_+$, Problem (\ref{approx_pb})
admits a unique family of solutions having the following structure, see also the former works \cite{GEORGIOU_LINDQUIST_KULLBACKLEIBLER}, \cite{PAVON_FERRANTE_ONGEORGIOULINDQUIST} and \cite{Ferrante-R-T-localconv},
\eq \Ac_{\hat \Theta,\nu}=\frac{\Psi}{(1+\nu^{-1}G^*\hat\Theta G)^\nu}\eeq
where $\hat\Theta \in\Qc_n$ is computed solving the dual problem
\eq \label{dual_alpha}\min\{J(\Theta),\; \Theta\in\Qc_n \hbox{ s.t. } 1+\nu^{-1}G^* \Theta G> 0\}\eeq
with
\eqn && J(\Theta)=\nn\\ && \left\{
                        \begin{array}{ll}
                          \int \Psi\log(1+G^*\Theta G)^{-1}+\tr(\hat\Sigma \Theta), & \nu=1 \\
                         \frac{\nu}{\nu-1}\int\Psi(1+\nu^{-1}G^*\Theta G)^{1-\nu}+\tr(\hat\Sigma \Theta) , & \nu>1.
                        \end{array}
                      \right. \eeqn
The multivariate case $m>1$, however, cannot be tackled with such divergence family. On the other hand, for the case $m=1$ and $\nu=2$, (\ref{def_alpha_divergence}) is the {\em Hellinger} distance and it can be extended to the case $m>1$ in such a way
Problem (\ref{approx_pb}) can be addressed, see \cite{FERRANTE_PAVON_RAMPONI_HELLINGERVSKULLBACK} and
\cite{RAMPONI_FERRANTE_PAVON_GLOBALLYCONVERGENT}.

\subsection{Solution with the Beta Divergence Family}
Consider the Beta Divergence Family \citep{BETA}
\eqn && \label{beta_div} \Sc_B^{(\beta)}(\Phi\|\Psi)=\tr\int [ \frac{1}{\beta(\beta-1)}\Phi^\beta-\frac{1}{\beta-1}\Phi\Psi^{\beta-1}\nn\\ &&\hspace{1cm}+\frac{1}{\beta}\Psi^\beta],\;\; \beta\in \Rs\setminus \{0,1\} .\eeqn
For $\beta\rightarrow 0$ and $\beta\rightarrow 1$ it is defined  by continuity
\eqn &&\lim_{\beta\rightarrow 0}\Sc_B^{(\beta)}(\Phi\|\Psi)=\Sc_{IS}(\Phi\|\Psi) \nn\\
&&\lim_{\beta\rightarrow 1}\Sc_B^{(\beta)}(\Phi\|\Psi)=\Sc_{KL}(\Phi\|\Psi)\eeqn
Here, $\Sc_{IS}$ denotes the {\em Itakura-Saito} distance
\eq \Sc_{IS}(\Phi\|\Psi)=\tr\int [\log\Psi-\log\Phi+\Phi\Psi^{-1}-I] .\eeq
In \cite{BETA} and \cite{FERRANTE_PAVON_MASIERO_SPECTRALRELATIVEENTROPY}, it has been shown that Problem (\ref{approx_pb}) with the parametrized divergence family $\Sc_B^{(1-\nu^{-1})}(\Phi\|\Psi)$, $\nu\in\Ns_+$,
admits a unique family of solutions of the form \eq \Bc_{\hat \Theta,\nu}=(\Psi^{-\nu^{-1}}+\nu^{-1}G^*\hat \Theta G)^{-\nu}\eeq
where $\hat\Theta\in\Qc_n$ is computed through the dual problem:
\eq \label{dual_beta}\min\{J(\Theta),\; \Theta\in\Qc_n \hbox{ s.t. } \Psi^{-\nu^{-1}}+\nu^{-1}G^* \Theta G\succ 0\}\eeq
with
\eqn && J(\Theta)=\nn\\ && \left\{
                        \begin{array}{ll}
                          \tr \int \log(\Psi^{-1}+G^*\Theta G)^{-1}+\tr(\hat\Sigma \Theta), & \nu=1 \\
                         \frac{\nu}{\nu-1}\tr\int (\Psi^{-\nu^{-1}}+\nu^{-1}G^*\Theta G)^{1-\nu}+\tr(\hat\Sigma \Theta) , & \nu>1.
                        \end{array}
                      \right. \nn \eeqn
The above solution also holds for the multivariate case, i.e. $m>1$, however, it requires the additional assumption that $\Psi^{\nu^{-1}}$ has bounded {\em McMillan} degree. Finally, it is worth noting that the limit case $\beta\rightarrow 1$ has been addressed in \cite{GEORGIOU_RELATIVEENTROPY} and it represents the first THREE-like method for the multivariate case.

\subsection{Solution with the Tau Divergence Family} \label{sec_sol_tau_div}
In \cite{BETAPRED}, it has been proposed the Tau divergence family
\eqn \label{tau_divergence} &&  \Sc_T^{(\tau)}(\Phi\|\Psi)=\tr\int[\frac{1}{\tau(\tau-1)}(W_\Psi^{-1}\Psi W_\Psi^{-*})^\tau\nn\\ && \hspace{1cm}-\frac{1}{\tau-1}\Phi\Psi^{-1}+\frac{1}{\tau}I],\hspace{0.5cm}\tau\in\Rs\setminus\{0,1\}\eeqn
where $W_\Psi$ is a left squared spectral factor of $\Psi$, that is $\Psi=W_\Psi W_\Psi^*$.
Moreover, for $\tau\rightarrow 0$ and $\tau\rightarrow 1$ we have
\eqn &&\lim_{\tau\rightarrow 0}\Sc_T^{(\tau)}(\Phi\|\Psi)=\Sc_{IS}(\Phi\|\Psi) \nn\\
&&\lim_{\tau\rightarrow 1}\Sc_T^{(\tau)}(\Phi\|\Psi)=\Sc_{KL}(W_\Psi^{-1}\Phi W_\Psi^{-*}\|I) .\eeqn
Problem (\ref{approx_pb}) with $\Sc_T^{(1-\nu^{-1})}(\Phi\|\Psi)$, $\nu\in\Ns_+$, admits a unique family of solutions of the form
\eq \label{opt_form_phi_T} \Tc_{\hat \Theta,\nu}=W_\Psi(I+\nu^{-1} W_\Psi^* G^*\hat \Theta G W_\Psi)^{-\nu} W_\Psi^*\eeq
and $\hat \Theta  \in\Qc_n$ is given solving the dual problem
\eq \label{dual_tau}\min\{J(\Theta),\; \Theta\in\Qc_n \hbox{ s.t. } I+\nu^{-1} W_\Psi^*G^* \Theta G W_\Psi\succ 0\}\eeq
with
\eqn && J(\Theta)=\nn\\ && \left\{
                        \begin{array}{ll}
                          \tr \int \log(\Psi^{-1}+G^*\Theta G)^{-1}+\tr(\hat\Sigma \Theta), & \nu=1 \\
                         \frac{\nu}{\nu-1}\tr\int(I+W_\Psi^* G^*\Theta G W_\Psi)^{1-\nu}+\tr(\hat\Sigma\Theta) , & \nu>1.
                        \end{array}
                      \right.  \eeqn Note that, the above solution holds for $m\geq 1$ under the mild assumption that $\Psi$ has bounded {\em McMillan} degree. Moreover, $ \Tc_{\hat \Theta,1}=\Bc_{\hat \Theta,1}$.
Finally, it is worth noting that \eq \Sc_T^{(1-\nu^{-1})}(\Phi\|\Psi)=\Sc_B^{(1-\nu^{-1})}(W_\Psi^{-1}\Phi W_\Psi^{-*}\| I). \eeq Here,
$W_\Psi^{-1}\Phi W_\Psi^{-*}$ is the spectral density of the normalized
prediction error process $\tilde \varepsilon=\{\tilde\varepsilon (t)\}_{t\in\Zs}$ 
where the actual model has spectral density $\Phi$ and the one-step ahead predictor is based on the prior $\Psi$. Accordingly, $\Sc_T^{(1-\nu^{-1})}(\Phi\|\Psi)$ represents a way to measure the mismatch between $\tilde \varepsilon$ and WGN. Therefore, (\ref{opt_form_phi_T}) is also the spectral density matching $\hat \Sigma$ and minimizing the prediction error $\tilde\varepsilon$. Problem (\ref{approx_pb}), however, cannot be reformulated as PEM. Indeed, in the latter the prediction error is optimized designing the one-step ahead predictor, see Section \ref{sec_PEM},  rather than the shaping filter of the process.

\section{Weighted Beta Divergence Families}\label{sec_weighted_beta}

Before to introduce our interpretation of the dual problem of the THREE-like approaches of Section \ref{sec_review_three}, we need to define two different types of Beta divergence weighted according to a weight function $Q\in\Sc_m^+$. Indeed, we will see in Section \ref{sec_dual_interpret} that the minimization of the dual function $J(\Theta)$ is equivalent to the minimization of a suitable weighted Beta divergence family. The latter measures the closeness between $\Omega$ and the THREE-like solution. 

\subsection{First type}
We can define the weighted Beta divergence as follows
\eqn && \Sc_{B1,Q}^{(\beta)}(\Phi\|\Psi)=\Sc_{B}^{(\beta)}(W_Q^*\Phi W_Q\|W_Q^*\Psi W_Q),\nn\\ && \hspace{5.4cm} \beta\in \Rs\setminus \{0,1\} \eeqn
where $W_Q$ is a left squared spectral factor of $Q$, i.e. $Q=W_Q W_Q^*$. Clearly, by choosing $Q=I$ we obtain the usual Beta divergence family defined in (\ref{beta_div}).
\begin{prop} For $\beta\in\Rs\setminus\{0,1\}$ and $Q\in\Sc_m^+$ fixed, $\Sc_{B1,Q}^{(\beta)}(\Phi\|\Phi)$ is a divergence index. Moreover, for $\beta\rightarrow 0$ and $\beta\rightarrow 1$ it can be extended by continuity
\eqn && \lim_{\beta \rightarrow 0} \Sc_{B1,Q}^{(\beta)}(\Phi\|\Psi)=\Sc_{IS}(\Phi\|\Psi)\nn\\
&& \lim_{\beta \rightarrow 1} \Sc_{B1,Q}^{(\beta)}(\Phi\|\Psi)=\Sc_{KL1,Q}(\Phi\|\Psi)\eeqn
where $\Sc_{KL1,Q}$ is the weighted Kullback-Leibler divergence
\eq \Sc_{KL1,Q}(\Phi\|\Psi)=\Sc_{KL}(W_Q^*\Phi W_{Q}\|W_{Q}^*\Psi W_Q).\eeq
\end{prop}
{\em Proof. } The statement can be proved by using the same lines of Proposition 2.1 in \cite{BETAPRED}. The unique difference regards the limit $\beta\rightarrow 0$:
\eqn && \lim_{\beta \rightarrow 0} \Sc_{B1,Q}^{(\beta)}(\Phi\|\Psi)=\lim_{\beta \rightarrow 0} \Sc_{B}^{(\beta)}(W_Q^*\Phi W_Q\| W_Q^*\Psi W_Q)\nn\\
&&\hspace{0.5cm} =\Sc_{IS}(W_Q^*\Phi W_Q\| W_Q^*\Psi W_Q)=\Sc_{IS}(\Phi\|\Psi)\eeqn
 where we exploited the property, see \cite{DISTANCES_JIANG_2012},
 \eq \label{proprieta_IS}\Sc_{IS}(\Phi_1\|\Phi_2)=\Sc_{IS}(W_{\Phi_2}^{-1}\Phi_1 W_{\Phi_2}^{-*}\| I)\eeq
 with $\Phi_1,\Phi_2\in\Sc_m^+$ and $\Phi_2=W_{\Phi_2}W_{\Phi_2}^*$.\qed\\
In view of (\ref{proprieta_IS}), it is worth noting that $\Sc_{IS}(\Phi\|\Psi)$ represents a way to measure the mismatch between the normalized prediction error 
$\tilde \varepsilon$ defined in Section \ref{sec_sol_tau_div} and WGN.
Finally, choosing $Q=\Psi^{-1}$ we obtain the Tau divergence family 
\eq \Sc_{B1,\Psi^{-1}}^{(\beta)}(\Phi\|\Psi)=\Sc_{T}^{(\beta)}(\Phi\|\Psi).\eeq

\subsection{Second type}
Another way to define the weighted Beta divergence follows
\eqn && \Sc_{B2,Q}^{(\beta)}(\Phi\|\Psi)=\tr\int Q[ \frac{1}{\beta(\beta-1)}\Phi^\beta-\frac{1}{\beta-1}\Phi\Psi^{\beta-1}\nn\\ &&\hspace{1cm}+\frac{1}{\beta}\Psi^\beta],\;\; \beta\in \Rs\setminus \{0,1\} .\eeqn Also in this case, setting $Q=I$ we obtain (\ref{beta_div}). 
\begin{prop} For $\beta\in\Rs\setminus\{0,1\}$ and $Q\in\Sc_m^+(\Ts)$ fixed, $\Sc_{B2,Q}^{(\beta)}(\Phi\|\Phi)$ is a divergence index. Moreover, for $\beta\rightarrow 0$ and $\beta\rightarrow 1$ it can be extended by continuity
\eqn && \lim_{\beta \rightarrow 0} \Sc_{B2,Q}^{(\beta)}(\Phi\|\Psi)=\Sc_{IS,Q}(\Psi\|\Phi)\nn\\
&& \lim_{\beta \rightarrow 1} \Sc_{B2,Q}^{(\beta)}(\Phi\|\Psi)=\Sc_{KL2,Q}(\Phi\|\Psi)\eeqn
where $S_{IS,Q}$ and $\Sc_{KL2,Q}$ are the weighted Itakura-Saito distance and the weighted Kullback-Leibler divergence, respectively,
\eq \Sc_{IS,Q}(\Phi\|\Psi)=\tr\int Q[\log\Psi-\log\Phi+\Phi\Psi^{-1}-I]\eeq
\eq \Sc_{KL2,Q}(\Phi\|\Psi)=\tr\int Q[\Phi(\log\Phi-\log\Psi)-\Phi+\Psi]\eeq
\end{prop}
{\em Proof.} Let $\beta\in\Rs\setminus\{0,1\}$. It is not difficult to show that $\Sc_{B2,Q}^{(\beta)}(\Phi\|\Psi)\succeq 0$ and $\Sc_{B2,Q}^{(\beta)}(\Psi\|\Psi)=0$. Assume that $\Sc_{B2,Q}^{(\beta)}(\Phi\|\Psi)=0$. Since $Q\in\Sc_m^+$,
there exists a positive constant $k$ such that $Q\succeq k I$. Moreover,
\eq k \Sc_{B}^{(\beta)}(\Phi\|\Psi)\preceq \Sc_{B2,Q}^{(\beta)}(\Phi\|\Psi)=0\eeq
accordingly $\Sc_{B}^{(\beta)}(\Phi\|\Psi)=0$ which implies $\Phi=\Psi$. We conclude that $\Sc_{B2,Q}^{(\beta)}$
is a divergence index. The same line can be exploited to show that $\Sc_{IS,Q}$ and $\Sc_{KL2,Q}$
are divergence indexes. Since it allowed to pass the limits $\beta\rightarrow 0$ and $\beta\rightarrow 1$ under the integral sign, see \cite{BETA}, then such limits can be easily proved by using Proposition 3.1 in \cite{BETA}.
\qed \\
Note that, for the scalar case $m=1$ we have
\eq \label{proprieta_itakura_saito_pesata}\Sc_{IS,Q}(\Phi\|\Psi)=\Sc_{IS,Q}(W_\Psi^{-1} \Phi W_\Psi^{-*}\| I)\eeq
therefore it represents a way to measure the mismatch between the normalized prediction error $\tilde \varepsilon$ defined in Section \ref{sec_sol_tau_div}
and WGN weighted according to the weight function $Q$. Finally, choosing $Q=\Psi^{1-\beta}$ we obtain the Alpha divergence family  \eq \Sc_{B2,\Psi^{1-\beta}}^{(\beta)}(\Phi\|\Psi)=\Sc_{A}^{(\beta)}(\Phi\|\Psi).\eeq

\section{The Dual Problem Interpretation}\label{sec_dual_interpret}
A THREE-like spectral estimator is solution to the spectrum approximation problem (\ref{approx_pb}). We now show its dual problem reveals this spectral estimator also solves another spectral estimation problem. More precisely, this estimator is the closest spectral density to the correlogram $\Omega$, according to a weighted Beta divergence family of Section \ref{sec_weighted_beta} , and belonging to a certain parametric class. We start by considering the dual problem (\ref{dual_tau}) with $\nu>1$. Taking into account (\ref{cond_correlogramma}), we obtain
\eqn && J(\Theta)=\tr\int\frac{\nu}{\nu-1}(I+\nu^{-1}W_{\Psi}^* G^*\Theta GW_{\Psi})^{1-\nu}\nn\\
&& \hspace{1.5cm}+\tr (\int G \Omega G^* \Theta)\nn\\
&& \hspace{1cm}=\tr\int[\frac{\nu}{\nu-1}(W_{\Psi}^{-1} \Tc_{\Theta,\nu} W_{\Psi}^{-*})^{1-\nu^{-1}}\nn\\
&& \hspace{1.5cm}+\nu \Omega W_{\Psi}^{-*}(\nu^{-1}W_{\Psi}^{*} G^* \Theta G W_{\Psi}) W_{\Psi}^{-1}]\nn\\\eeqn
where $\Tc_{\Theta,\nu}$ has been defined in (\ref{opt_form_phi_T}). Since the term $\tr\int \frac{\nu^2}{1-\nu}(W_\Psi^{-1}\Omega W_\Psi^{-*})^{1-\nu^{-1}}$ plays no role in the optimization with respect to $\Theta$, we can add it to $J$:
\eqn  && J(\Theta)=\tr\int[\frac{\nu}{\nu-1}(W_{\Psi}^{-1}\Tc_{\Theta,\nu} W_{\Psi}^{-*})^{1-\nu^{-1}}\nn\\
&& \hspace{1.5cm}+\nu \Omega W_{\Psi}^{-*}(I+\nu^{-1}W_{\Psi}^{*} G^* \Theta G W_{\Psi}) W_{\Psi}^{-1}\nn\\&& \hspace{1.5cm}+\frac{\nu^2}{1-\nu}(W_\Psi^{-1} \Omega W_\Psi^{-*})^{1-\nu^{-1}}]\nn\\
&& \hspace{1cm} =\tr\int[\frac{\nu}{\nu-1}(W_{\Psi}^{-1}\Tc_{\Theta,\nu} W_{\Psi}^{-*})^{1-\nu^{-1}}\nn\\
&& \hspace{1.5cm}+\nu W_{\Psi}^{-1}\Omega W_{\Psi}^{-*}(W_{\Psi}^{-1} \Tc_{\Theta,\nu} W_{\Psi}^{-*})^{-\nu^{-1}}\nn\\ &&\hspace{1.5cm} +\frac{\nu^2}{1-\nu}(W_\Psi^{-1}\Omega W_\Psi^{-*})^{1-\nu^{-1}}]\nn\\
&& \hspace{1cm} =\Sc_{B}^{(1-\nu^{-1})}(W_{\Psi}^{-1} \Omega W_{\Psi}^{-*}\|W_{\Psi}^{-1} \Tc_{\Theta,\nu}W_{\Psi}^{-*})\nn\\
&& \hspace{1cm} =\Sc_{B1,\Psi^{-1}}^{(1-\nu^{-1})}(\Omega \| \Tc_{\Theta,\nu}).\eeqn
Consider the parametric class of spectral densities
\eq \Mc_{T}=\{\Tc_{\Theta,\nu},\;\; \Theta\in \Dc_{T}\}\eeq
where the parameter matrix $\Theta\in\Dc_{T}=\{\Theta\in\Qc_n \hbox{ s.t. } I+\nu^{-1} W_\Psi^{*}G^*\Theta GW_\Psi\succ 0 \}$. Therefore, the dual problem (\ref{dual_tau}) is equivalent to
\eq \label{interpretazione_duale}\hat\Theta =\mathrm{argmin}_{\Theta\in\Dc_{T}} \Sc_{B1,\Psi^{-1}}^{(1-\nu^{-1})}(\Omega \|\Tc_{\Theta,\nu}).\eeq
Similarly, it can be proved that (\ref{interpretazione_duale}) also holds for the case $\nu=1$.
This interpretation of the dual problem allows to understand
$\Tc_{\hat \Theta,\nu}$ optimal in terms of a new spectral estimation problem.
\begin{prop}\label{prop_tau} $\Tc_{\hat \Theta,\nu}$ is the closest spectral density to the correlogram $\Omega$, according to $\Sc_{B1,\Psi^{-1}}^{(1-\nu^{-1})}(\Omega \|\Tc_{\Theta,\nu})$, and belonging to the parametric class $\Mc_{T}$.\end{prop}
Similar results can be derived for the estimators $\Ac_{\hat \Theta,\nu}$ and $\Bc_{\hat \Theta,\nu}$.
\begin{prop} \label{prop_alpha}$\Ac_{\hat \Theta,\nu}$  is the closest spectral density to the correlogram $\Omega$, according to $\Sc_{B2,\Psi^{\nu^{-1}}}^{(1-\nu^{-1})}(\Omega\| \Ac_{\Theta,\nu})$, and belonging to the parametric class \eq \Mc_{A}=\{\Ac_{\Theta,\nu},\;\; \Theta\in\Dc_{A})\}\eeq with $\Dc_{A}=\{\Theta\in\Qc_n \hbox{ s.t. } 1+\nu^{-1} G^*\Theta G>0\}$. \end{prop}
\begin{prop}\label{prop_beta} $\Bc_{\hat \Theta,\nu}$  is the closest spectral density to the correlogram $\Omega$, according to $\Sc_{B}^{(1-\nu^{-1})}(\Omega\|\Bc_{\Theta,\nu})$, and belonging to the parametric class \eq \Mc_{B}=\{\Bc_{\Theta,\nu},\;\; \Theta\in\Dc_{B}\}\eeq with $\Dc_{B}=\{\Theta\in\Qc_n \hbox{ s.t. } \Psi^{-\nu^{-1}}+\nu^{-1} G^*\Theta G\succ 0\}$. \end{prop}

We conclude that the solution to the spectrum approximation problem (\ref{approx_pb}) can be seen as the solution of a parametric spectral estimation problem wherein the best estimate is the closest one to the correlogram, according to an appropriate divergence index, and belonging to a suitable parametric class. Moreover, the {\em a priori} spectral density $\Psi$ always belongs to the parametric class.
Finally, it is worth noting that the dual problem is always characterized by Beta-like divergence families. Indeed, such divergence indexes are the unique to contain a linear term in $\Omega$ corresponding to the term $\tr(\Theta\hat\Sigma)$ in the dual function.

\section{The connection between THREE and PEM}\label{sec_PEM}
Consider the case $m=1$ and $\Psi=1$. Then, 
it is not difficult to see the spectral estimator $\Ac_{\hat \Theta,1}$ (which coincides with $\Bc_{\hat \Theta,1}$ and $\Tc_{\hat \Theta,1}$)
minimizes the divergence index $\Sc_{IS}(\Omega\| \Ac_{\hat \Theta,1})$. The latter is the relative entropy rate between two stationary Gaussian processes 
having spectral density $\Omega$ and $\Ac_{\hat \Theta,1}$, respectively \citep{COVER_THOMAS}. In 
\cite{LINDQUIST_PREDICTION_ERROR_2007} and \cite{BLOMQVIST_WAHLBERG_2007}
it has been shown that the minimization of the relative entropy rate is equivalent to the 
prediction error identification method (PEM). Therefore, $\Ac_{\hat \Theta,1}$ is also solution to PEM.
We now show this result can be extended to the case $\Psi\neq 1$ for $\Ac_{\hat \Theta,1}$ and to the case $m\geq 1$, $\Psi\neq I$ for $\Tc_{\hat \Theta,1}(=\Bc_{\hat \Theta,1})$.

 First, we review the PEM approach. Let $y=\{y(t)\}_{t\in\Zs}$ be a $\Rs^m$-valued, zero-mean, purely nondeterministic, full rank, stationary, Gaussian stochastic process having model
\eq \label{y_in_model_class}\Pc_\Theta \;: \; y(t)=\sum_{k=0}^\infty F_{\Theta,k}e(t-k)\eeq
where $e:=\{e(t)\}_{t\in\Zs}$ is the normalized innovation process, i.e. WGN,
and $\{F_{\Theta,k}\}_{k\in\Ns}$, $F_{\Theta,k}\in\Rs^{m\times m}$, is the impulse response of the shaping filter. Moreover, $\Pc_\Theta$ belongs to the class of models
\eq \Mc:=\{\Pc_\Theta \;|\; \Theta\in \Dc\}\eeq
wherein each model is parametrized using the parameter vector (or possibly matrix) $\Theta\in D\subset \Rs^d$. Let $\varepsilon_\Theta=\{\varepsilon_{\Theta}(t)\}_{t\in\Zs}$ be the normalized prediction error of $y$ where the 
one-step ahead predictor is based on model $\Pc_{\Theta}$. If $\Pc_\Theta$ is the true model for $y$, then $\varepsilon_{\Theta}=e$, that is $\varepsilon_\Theta$ is WGN. On the contrary, the closer $\varepsilon_{\Theta}$ is to be WGN, the better $\Pc_\Theta$ describes $y$.
Consider now the situation that a finite length sequence $\mathrm{y}:=\{\mathrm{y}(t)\}_{t=1}^N$  extracted form a realization of $y$ is given. We consider the problem to select an appropriate value $\hat \Theta  \in \Dc$ of the parameter vector, and therefore an appropriate model $\Pc_{\hat \Theta}\in\Mc$,  by using the information in $\mathrm{y}$. The prediction error identification method judges the performance in respect to the prediction error $\varepsilon_{\Theta}$ of each model $\Pc_\Theta\in\Mc$ and then selects
as $\hat \Theta$ the one with the best performance. More precisely, we consider the cost function $V(\Theta,\mathrm{y})$
which is a scalar-valued positive function of $\varepsilon_\Theta$.
Therefore, $\hat \Theta$ is obtained solving the following optimization problem
\eq \label{PEM_problem}\hat\Theta =\mathrm{argmin}_{\Theta\in \Dc} V(\Theta,\mathrm{y}). \eeq
If we choose \eq V(\Theta,\mathrm{y})=\frac{1}{N}\sum_{t=1}^N \|\varepsilon_{\Theta}(t)\|^2\eeq
we obtain the standard PEM \citep{LJUNG_IDENTIFICATIONTHEORYUSER}.
Let $L_{\Theta}$ be the {\em Fourier} transform of the sequence $\{ F_{\Theta,k}\}_{k\in\Ns}$. Then,
$\Phi_\Theta=L_\Theta L_\Theta^*$ is the spectral density of $y$ and is equivalent to $\Pc_\Theta$. Starting from this observation we
show that the models $\Ac_{\hat \Theta,1}$  and $\Tc_{\hat \Theta,1}$ ($=\Bc_{\hat \Theta,1}$)
can be seen as solution to (\ref{PEM_problem}).
\subsection{Solution $\Ac_{\hat \Theta,1}$}
By Proposition \ref{prop_alpha}, we known that
$\hat\Theta$ is given by the minimization of $\Sc_{IS,\Psi}(\Omega \|\Ac_{ \Theta,1})$ with $\Theta \in \Dc_A$. In view of (\ref{proprieta_itakura_saito_pesata}), we have
$\Sc_{IS,\Psi}(\Omega \|\Ac_{\Theta,1})=\Sc_{IS,\Psi}(L_\Theta^{-1}\Omega L_\Theta^{-*}\|I)$
where $L_\Theta$ is such that $\Ac_{\Theta,1}=L_\Theta L_\Theta^*$.
Let $\varepsilon_\Theta=\{\varepsilon_\Theta(t)\}_{t\in\Zs}$ be the normalized prediction error process 
where the actual process has spectral density $ \Omega$ and
the one-step ahead predictor is  based on the model
$\Ac_{\Theta,1}$. It is not difficult to see that
$\Lambda_\Theta=L_\Theta^{-1}\Omega L_\Theta^{-*}$
denotes the spectral density of $\varepsilon_\Theta$. Accordingly,
by choosing
 \eq V(\Theta,\mathrm{y})=\Sc_{IS,\Psi}(\Lambda_\Theta \| I)\eeq
we obtain the PEM problem
\eq \hat\Theta=\mathrm{argmin}_{\Theta\in \Dc_A} V(\Theta,\mathrm{y}).\eeq

\subsection{Solution $\Tc_{\hat \Theta,1}$}
By Proposition \ref{prop_tau}, we have that
$\hat\Theta $ is given by the minimization of $\Sc_{IS}(\Omega \| \Tc_{\Theta,1})$ with $\Theta\in \Dc_{T}$.
Moreover, \eqn &&\Sc_{IS}(\Omega \|\Tc_{\Theta,1})=\Sc_{IS}(L_\Theta^{-1}\Omega L_\Theta^{-*}\|I)\eeqn
where $L_\Theta$ is such that $\Tc_{\Theta,1}=L_\Theta L_\Theta^*$.
Therefore, similarly to the previous case, we have
\eq \label{PEM_itakura_saito} \hat\Theta  =\mathrm{argmin}_{\Theta\in \Dc_{T}} V(\Theta,\mathrm{y})\eeq
 where
 \eq V(\Theta,\mathrm{y})=\Sc_{IS}(\Lambda_\Theta \| I).\eeq
 Here, $\Lambda_\Theta$ is the spectral density of the normalized prediction error process where the 
 actual process has spectral density $\Omega$ and the one-step ahead predictor is based on $\Tc_{\Theta,1}$.
Note that, in (\ref{PEM_itakura_saito})  we can replace $\Dc_T$ with $\Dc_T\cap \Vc$ where $\Vc$ is a vector subspace of $\Qc_n$. Accordingly, $\Dc_T\cap \Vc \subset \Dc_T$. At this point, recall that a $\Cs^{m\times m}$-valued analytic matrix function is sparse if many of its entries are null functions, and is low-rank if its pointwise rank (constant almost everywhere) is low almost everywhere. By choosing properly $G$, $\Psi$ and $\Vc$ the parametric class of models \eq \Mc_{T,\Vc}=\{\Tc_{\Theta,1},\;\;\Theta\in \Dc_T\cap \Vc\}\eeq
only contains spectral densities whose inverse is sparse \citep{ARMA_AVVENTI_LINDQUIST_WAHLBERG_2013} or sparse plus low rank \citep{LATENTG}. Such parametric classes are important in graphical modeling \citep{LAURITZEN_1996} where the process $y$ is ``attached'' to a graph: each node corresponds to a variable in $y$ and there is a direct link between two variables if and only if are conditional dependent given the remaining variables. We conclude that also the solutions presented in \cite{ARMA_AVVENTI_LINDQUIST_WAHLBERG_2013} and in \cite{LATENTG} admit a PEM interpretation similar to (\ref{PEM_itakura_saito}).

\section{Conclusions}
In this paper, we have presented an interpretation of the dual problem arising from the THREE-like methods. From this interpretation it turns out that the solution to a THREE-like problem is also the closest spectral density to the correlogram over a certain parametric class. Moreover, two particular solutions can be seen also as solution to PEM.

\end{document}